\documentclass[12pt]{amsart}
\usepackage{amsmath}
\usepackage{amsfonts,amssymb,amsthm}
\usepackage{graphicx}

\newcommand{\beq}{\begin{equation}}
\newcommand{\eeq}{\end{equation}}
\newcommand{\bbar}{\begin{eqnarray}}
\newcommand{\eear}{\end{eqnarray}}

\newcommand{\thm}[2]{\begin{#1} #2 \end{#1}}

\newtheorem{theorem}{Theorem}[section]

\newtheorem{lemma}[theorem]{Lemma}

\begin{document}

\title{Mathematics of Learning}

\author{Natalia L. Komarova}
\address{Institute For Advanced Study, Einstein Drive, Princeton, NJ
$08540$}
\address{Department of Applied Mathematics, University of
Leeds, Leeds LS2 9JT, UK}

\email{natalia@ias.edu}
\author{Igor Rivin}
\address{Mathematics department, University of Manchester,
Oxford Road, Manchester M13 9PL, UK}
\address{Mathematics Department, Temple University,
Philadelphia, PA 19122}
\address{Mathematics Department, Princeton University, Princeton,
NJ 08544}
\email{irivin@math.princeton.edu} \thanks{N.Komarova gratefully
acknowledges support from the Packard Foundation, the Leon Levy and
Shelby White Initiatives Fund, the Florence Gould Foundation, the
Ambrose Monell Foundation, the Alfred P. Sloan Foundation and the
NSF. I. Rivin would like to thank the EPSRC and the NSF for support}

\subjclass{60E07, 60F15, 60J20, 91E40, 26C10} \keywords{harmonic 
mean, random polynomials, random matrices, learning theory, 
Markov processes, stable law, large deviations}
\begin{abstract}
We study the convergence properties of a pair of learning algorithms
(learning with and without memory). This leads us to study the
dominant eigenvalue of a class of random matrices. This turns out to
be related to the roots of the derivative of random polynomials
(generated by picking their roots uniformly at random in the interval
$[0, 1]$, although our results extend to other distributions). This,
in turn, requires the study of the statistical behavior of the
harmonic mean of random variables as above, which leads
us to delicate question of the rate of convergence to stable laws and
tail estimates for stable laws. The reader can find the proofs of most
of the results announced here in \cite{kr1}.
\end{abstract}
\maketitle

\renewcommand{\theitheorem}{\Alph{itheorem}}
The original motivation for the work in this paper was provided by the
first-named author's research in learning theory, specifically in
various models of language acquisition (see \cite{knn,nkn,kn}) and
more specifically yet by the analysis of the speed of convergence of
the \emph{memoryless learner algorithm}. Curiously, our methods also
result in a complete analysis of \emph{learning with full memory}, as
shown in some detail in section \ref{fullmem}. The setup is described
in section \ref{memoryless}, so here we will just recall the
essentials. There is a collection of concepts $R_1, \dots, R_n$ and
words which refer to these concepts, sometimes ambiguously. The teacher
generates a stream of words, referring to the concept $R_1$. This is
not known to the student, but he must learn by, at each steps,
guessing some concept $R_i$ and checking for consistency with the
teacher's input.  The memoryless learner algorithm consists of picking
a concept $R_i$ at random, and sticking by this choice, until it is
proven wrong.  At this point another concept is picked randomly, and
the procedure repeats. \emph{Learning with full memory} follows the
same general process with the important difference that once a concept
is rejected, the student never goes back to it\footnote{Another
important learning algorithm is the so-called \emph{batch
learner}. This is analysed completely in \cite{ri1}}. It is clear that once
the student hits on the right answer $R_1$, this will be his final
answer, so the question is then:
\begin{quote}
\emph{How quickly do the two methods converge to the truth?}
\end{quote}
Since the first method is \emph{memoryless}, as the name implies, it
is clear that the learning process is a \emph{Markov process}, and as
is well-known the convergence rate is determined by the gap between
the top (Perron-Frobenius) eigenvalue and the second largest
eigenvalue.  However, we are also interested in a kind of a
\emph{generic} behavior, so we assume that the sizes of overlaps
between concepts are \emph{random}, with some (sufficiently regular)
probability density function supported in $[0, 1]$, and that the
number of concepts is large. This makes the transition matrix random,
though of a certain restricted kind, as described in detail in section
\ref{memoryless}. The analysis of convergence speed then comes down to
a detailed analysis of the size of the second-largest eigenvalue and
also of the properties of the eigenspace decomposition. The analysis
for learning with full memory is quite different, but the results have
a very similar form. We summarize below:
\thm{theorem}
{Let $N_\Delta$ be the number of steps it takes for the student 
to have probability $1 - \Delta$ of learning the
concept. Then we have the following estimates for $N_\Delta$:
\begin{itemize}
\item
if the distribution of overlaps is \emph{uniform}, or more 
generally, the density function $f(1-x)$  at $0$ has the form 
$f(x) = c + O(x^\delta),$ $\delta, c > 0,$ then there exist positive
constants $C_1, C_2, C_1', C_2'$ such that 
$$\lim_{n \rightarrow \infty} \mathbf{P}\left(C_1 <
\frac{N_\Delta}{|\log \Delta|n \log n} < C_2\right) = 1$$
for
the memoryless algorithm and
$$\lim_{n \rightarrow \infty} \mathbf{P}\left(C'_1 <
\frac{N_\Delta}{(1- \Delta)^2 n \log n} < C'_2\right) = 1$$
when learning with full memory;
\item 
if the probability density function $f(1-x)$ is asymptotic to $c x^\beta
+ O(x^{\beta + \delta}), \quad \delta, \beta > 0$, as $x$ approaches
$0$, then for the two algorithms we have respectively 
$$\lim_{n \rightarrow \infty} \mathbf{P}\left(c_1 <
\frac{N_\Delta}{|\log \Delta|n} < c_2\right) = 1,$$
and
$$\lim_{n \rightarrow \infty} \mathbf{P}\left(c_1' <
\frac{N_\Delta}{(1- \Delta)^2 n } < c_2'\right) = 1$$
for some positive constants $c_1, c_2, c_1', c_2'$;
\item 
if the asymptotic behavior is as above, but $-1 < \beta < 0$, then
$$\lim_{x \rightarrow \infty}  \mathbf{P}\left(\frac{1}{x} < \frac{N_\Delta}{|\log
\Delta| n^{1/(1+\beta)}} < x\right) = 1$$
for the memoryless learning algorithm, and similarly
$$\lim_{x \rightarrow \infty}  \mathbf{P}\left(\frac{1}{x} < \frac{N_\Delta}{(1-\Delta)^2 n^{1/(1+\beta)}} < x\right) = 1$$
for learning with full memory.
\end{itemize}}
It should be said that our methods give quite precise estimates on the
constants in the asymptotic estimate, but the rate of convergence is
rather poor -- logarithmic -- so these precise bounds are of limited
practical importance.

\section{Eigenvalues and polynomials} In order to calculate the
convergence rate of the learning algorithm described above, we 
need to study the spectrum of a class of random matrices. The 
matrices have the following form:
\begin{equation}
\label{Tmatrix}
T_{ij} = 
  \begin{cases}
    a_i & i=j, \\
    \frac{(1-a_i)}{n-1} & \text{otherwise}, 
  \end{cases}
\end{equation}
where 
\begin{equation}
\label{valuesa}
a_1=1,\quad\quad 0\le a_i<1,\quad 2\le i\le n.
\end{equation}
Let $B=\frac{n-1}{n}(I-T)$, so that the eigenvalues of $T$,
$\lambda_i$, are related to the eigenvalues of $B$, $\mu_i$ by
$\lambda_i=1-\left[n/(n-1)\right]\mu_i$. We show the following amusing
\thm{lemma}{
\label{charpd} Let $p(x) = (x - x_1) \dots (x-x_n)$, where 
$x_i=1-a_i$. Then the characteristic polynomial $p_B$ of $B$ 
satisfies:
$$
p_B(x) = \frac{x}{n} \frac{dp(x)}{dx}.
$$
}
From lemma \ref{charpd}, the second largest eigenvalue of the matrix
$T$, $\lambda_*$, and the smallest root of $p'(x)$, which we denote as
$\mu_*$, are related as
\begin{equation}
\label{muandlamb}
\lambda_*=1-\frac{n}{n-1}\mu_*.
\end{equation}
Therefore, we need to study the distribution of the \emph{smallest}
root of $p^\prime(x)$, given that the smallest root of $p(x)$ is fixed
at $0$. Letting the roots of $p(x)$ be $0 = x_1 < x_2 < \dots <
x_{n}$, and letting 
\begin{equation}
\label{defharmean}
H(x_2, \dots, x_n)=\frac{(n-1)}{\sum_{i=2}^n1/x_i}
\end{equation}
be the \emph{harmonic mean} of the nontrivial roots of $p(x)$, we have
\thm{theorem}{
\label{ismall}
 The smallest root $\mu_*$ of $p^\prime(x)$
satisfies:
\begin{equation}
\label{Htomu}
\frac{1}{2}H(x_2, \dots, x_n) \leq (n-1) \mu_* \leq H(x_2, \dots, x_n).
\end{equation}
} 
We can see that the study of the distribution of $\mu_*$ entails the
study of the distribution of the asymptotic behavior of the harmonic
mean of a sample from a distribution on $[0, 1]$.

\section{Statistics of the harmonic mean.} 
The arithmetic, harmonic, and geometric means are examples of the
``conjugate means'', given by
\begin{equation*}
m_{\mathcal F}(x_1, \dots, x_n) = {\mathcal F}^{-1} \left(\frac{1}{n} \sum_{i=1}^n 
{\mathcal F}(x_i)\right),
\end{equation*}
where ${\mathcal F}(x) = x$ for the arithmetic mean, ${\mathcal F}(x)
= \log(x)$ for the geometric mean, and ${\mathcal F}(x) = 1/x$ for the
harmonic mean. The interesting situation is when ${\mathcal F}$ has a
singularity in the support of the distribution of $x$, and this case
seems to have been studied very little, if at all. Here we will devote
ourselves to the study of harmonic mean.

Given $x_1, \dots, x_n$ -- a sequence of independent, identically
distributed in $[0, 1]$ variables (with common probability density
function $f$), the nonlinear nature of the harmonic mean leads us to
consider first the random variable
\begin{equation}
\label{Xdef}
X_n = \frac{1}{n} \sum_{i=1}^n \frac{1}{x_i}. 
\end{equation}
Since the variables $1/x_i$ are easily seen to have infinite
expectation and variance, our prospects seem grim at first 
blush, but then we notice that the variable $1/x_i$ falls 
straight into the framework of the ``stable laws'' of L\'evy -- 
Khintchine, which is briefly presented below. 

\subsection{\label{stablaw}Stable limit laws}
Consider an infinite sequence of independent identically distributed
random variables $y_1, \dots, y_n, \dots$, with probability
distribution function $\mathfrak{F}$. Typical questions studied in
probability theory are the following.
\begin{quote}
Let $S_n = \sum_{j=1}^n y_j$. How is $S_n$ distributed? What can 
we say about the distribution of $S_n$ as $n \rightarrow \infty$?
\end{quote}
The best known example is one covered by the Central Limit 
Theorem: if $\mathfrak{F}$ has finite mean 
$\mu$ and variance $\sigma^2$, then $(S_n - n \mu)/(\sqrt{n} \sigma)$ 
converges in distribution to the \emph{normal distribution} 
(\cite{nor}). Similarly, we say that the variable 
$X$ \emph{belongs to the domain of attraction of a non-singular 
distribution $G$}, if there are constants $a_1, \dots, a_n, 
\dots$ and $b_1, \dots, b_n, \dots$ such that the sequence of 
variables $Y_k=a_k S_k - b_k$ converges in distribution to $G$. 
It was shown by L\'evy  and by Khintchine that having a domain of 
attraction constitutes severe restrictions on the distribution as 
well as the {\it norming sequences} $\{a_k\}$ and $\{b_k\}$. To 
wit, one can always pick $a_k = k^{- 1/\alpha}l(k), \quad 0 < \alpha 
\leq 2,$ where $l(k)$ is a slowly varying function (in the sense of
Karamata). In that case, $G$ is called a \emph{stable distribution   
of exponent $\alpha$}. If the variable $y$ belongs to the domain 
of a stable distribution of exponent $\alpha > 1$, then $y$ has 
an expectation $\mu$; just as in the case $\alpha = 2$, we can 
choose $b_k = k^{1-1/\alpha} \mu.$ When $\alpha < 1$, the 
variable $y$ has no mean, and it turns out that we can take $b_k 
\equiv 0$; for $\alpha = 1$, we can take $b_n = c \log n$, where 
$c$ is some constant depending on $\mathfrak{F}$. In particular, 
the normal distribution is a stable distribution of exponent $2$ 
(and is unique, up to scale and shift). This is one of the few 
cases where we have an explicit expression for the density of a 
stable distribution. The Fourier transforms of the densities are
explicitely known; the reader can find them in \cite[Chapter
XVII]{feller}. The stable distribution of a given exponent are 
parameterized by parameters $p, q, C$ defined below:
\begin{eqnarray}
\label{righttail} \
&&\lim_{x \rightarrow \infty} \frac{1-\mathfrak{F}(x)}{1-\mathfrak{F}(x)
+ \mathfrak{F}(-x)} = Cp,\\
\label{lefttail}
&&\lim_{x\rightarrow \infty}\frac{\mathfrak{F}(-x)}{1-\mathfrak{F}(x)} = Cq, 
\end{eqnarray}
and $p+q=1$. We will say that the stable law is \emph{unbalanced} if
$p=1$ or $q=1$ above. This will happen if the support of the variable
$y$ is positive -- this will be the only case we will consider in the
sequel.  Note that this \emph{does not} mean that the stable
distribution is supported away from $-\infty$, though that is true for
exponents smaller than $1$.
\subsection{\label{distharm}Limiting distribution of the harmonic mean}
Which particular stable law comes up in the study of the variable
$X_n$ in (\ref{Xdef}), depends on the distribution function $f(x)$. Let us
assume that 
$$f(x) \asymp c x^\beta,$$
as $x \rightarrow 0$ (for the uniform distribution $\beta = 0, \quad c
= 1$). (The notation $b \asymp a$ means that $a$ is asymptotically the
same as $b$, i.e.  there exist constants $c_1$, $c_2$, $d_1$, $d_2$,
so that $c_1 b + d_1 \leq a \leq c_2 b + d_2$.)  Then we have
\thm{theorem}{ \label{stableob} If $\beta = 0$, then let $Y_n 
= X_n - \log n.$ The variables $Y_n$ converge in distribution to 
the variable $\mathcal{Y}$ distributed in accordance to the unbalanced
stable law $G^{(\alpha)}$ with $\alpha=1$. If  
$\beta > 0$, then $X_n$ converges in distribution to 
$\delta(x-\mu)$, where $\mu = \mathbf{E}(1/x_i)$ (since the $x_i$ are
identically distributed the value of the index $i$ is not
relevant). If $-1 < \beta < 0$, then $n^{1-1/(1+\beta)} X_n$ converges
in distribution to a the variable $\mathcal{X}$ distributed in
accordance to the  stable law with exponent $\alpha=1+\beta$. }
\par\medskip\noindent \textbf{Remark.} In the case when the variables
$x_1, \dots, x_i$ have positive and continuous density at $0$, the
variables $X_n$ above converge to the Cauchy distribution (the
symmetric stable distribution of exponent $1$). This is the content of
exercise 7.6 in \cite{durrett}, though the (necessary) condition of
positivity of the density at $0$ is inadvertently omitted there.

The Theorem \ref{stableob} points us in the right direction, since it
allows us to guess the form of the following results ($H_n$ is the
harmonic mean of the variables):
\thm{theorem}
{
\label{harmlog}
Let $H_n = 1/X_n$ and $\beta = 0$. Then there exists a constant $\mathfrak{C}_1$ such that
\begin{equation*}
\lim_{n\rightarrow \infty} 
 \mathbf{E}(H_n\log n ) = \mathfrak{C}_1.
\end{equation*}
}

\begin{theorem}
\label{harmmu}
Suppose $\beta > 0$, let $y = 1/x$, and let $\mu$ be the mean of the
variable $y.$ Then  $\mathbf\lim_{n\rightarrow
\infty} \mathbf{E}(\mu H_n) = 1$.
\end{theorem}
Finally,
\begin{theorem}
\label{expneg}
Suppose $\beta < 0.$ Then there exists a constant $\mathfrak{C}_2$ such
that $\mathbf{E}(H_n/n^{1-1/{(1+\beta)}}) = \mathfrak{C}_2.$
\end{theorem}
We also have the following laws of large numbers:

\thm{theorem}{\textbf{Laws of large numbers for harmonic mean.}
\label{iweak} Let $\beta = 0$ and let $a  
> 0$. Then 
\begin{equation*}
\lim_{n\rightarrow \infty} \mathbf{P}( \vert H_n \log n -\mathfrak{C}_1\vert 
> a) = 0.
\end{equation*}
If $\beta > 0$, and $\mu$ is as in the statement of Theorem
\ref{harmmu},  then 

\begin{equation*}
\lim_{n\rightarrow \infty} \mathbf{P}(\vert H_n 
-\frac{1}{\mu}\vert) > a) = 0.
\end{equation*}}

The proofs of the above results use a variety of estimates; the reader
is referred to \cite{kr1}.  In addition to the laws of large numbers,
we also have the following limiting distribution results:

\begin{theorem}
\label{hlimlaw} For $\alpha=1$, the random variable $\log n(H_n\log n 
- {\mathfrak C}_1)$ converges to $1-G(-x/{\mathfrak C}_1^2),$ where $G$
is the limiting distribution (of exponent $\alpha=1$) of variables
$Y_n = X_n - c\log n$ and ${\mathfrak C}_1=1/c$.
\end{theorem}

\begin{theorem}
\label{hlimlaw2}
For $\alpha>1$, the random variable $n^{1-1/\alpha} (H_n -
\frac{1}{{\mathcal E}})$ converges in distribution to the variable
$\mathcal{H}$ with the distribution function $1-G(-x {\mathcal
E}^2)$, where $G$ is the unbalanced stable law of exponent $\alpha$.
\end{theorem}

\begin{theorem}
\label{hlimlaw3}
For $0<\alpha<1$, the random variable $H_n/n^{1-1/\alpha}$ converges in
distribution to the variable $\mathcal{H}$, with the distribution function
$1-G(1/x)$, where $G$ is the distribution function of the unbalanced
stable law of exponent $\alpha.$
\end{theorem}

\bibliographystyle{alpha}

\section{A pair of learning algorithms}
\subsection{The memoryless learner algorithm}
\label{memoryless}
Suppose there are $n$ intersecting sets, $R_1,\ldots, 
R_n$, and $n$ probability measures, $\nu_1,\ldots, \nu_n$, each 
defined on its set (so that $\nu_i(R_i)=1$). The 
\textit{similarity matrix} $A$ is given by $a_{ij}=\nu_i(R_j)$. It 
follows that $0\le a_{ij}\le 1$ and $a_{ii}=1$ for all $i$ and 
$j$.

Let us consider a typical problem of learning theory. A teacher
generates a sequence of points which belong to one of these sets, say
to set $R_1$. The total length of the sequence is $N$. The learner's
task is to guess what set is the teacher's set after receiving $N$
points. For simplicity we assume here that $a_{ij}<1$ for $i\ne j$,
which means that no set is a subset of another set. Many different
algorithms are available to the learner, one given by the so-called
{\it memoryless learner} algorithm \cite{niy}, a favorite with
learning theorists. It works in the following way. The learner starts
by (randomly) choosing one of the $n$ sets as an initial state. Then
$N$ sample points are received from the teacher. For each sampling,
the learner checks if the point belongs to its current set.  If it
does, no action is taken; otherwise, the learner randomly picks a
different set. The initial probability distribution of the learner is
uniform: ${\bf p}^{(0)}=(1/n,\ldots,1/n)^T$, i.e. each of the sets has
the same chance to be picked at the initial moment. The discrete time
evolution of the vector ${\bf p}^{(t)}$ is a Markov process with
transition matrix $T$, which depends on the similarity matrix,
$A$. The transition matrix is given by Eqs. (\ref{Tmatrix}),
(\ref{valuesa}) with $a_i=\nu_1(R_i)$.

After $N$ samplings, the probability of learning the correct set is
given by $\label{Qij} Q_{11}=[({\bf p}^{(0)})^T\,T^N]_1.$ It is clear
that the convergence rate of the memoryless algorithm can be
determined if we study properties of the matrix $T$. We are interested
in the rate of convergence as a function of $n$, the number of
possible sets.

We define the convergence rate of the method as the difference
$1-Q_{11}$. In order to evaluate the convergence rate of the
memoryless learner algorithm, let us represent the matrix $T$ as
$T=V\Lambda W,$ where the diagonal matrix $\Lambda$ consists of the
eigenvalues of $T$, which we call $\lambda_i$, $1\le i\le n$; the
columns of the matrix $V$ are the right eigenvectors of $T$, ${\bf
v}_i$, and the rows of the matrix $W$ are the left eigenvectors of $T$,
${\bf w}_i$, normalized to satisfy $<{\bf w}_i, {\bf
v}_j>=\delta_{ij}$ (so that $VW=WV=I$). The eigenvalues of $T$ satisfy
$|\lambda_i|\le 1$. We have
\begin{equation*}
T^N=V\Lambda^NW.
\end{equation*}
Let us arrange the eigenvalues so that $\lambda_1=1$ and
$\lambda_2\equiv \lambda_*$ is the second largest eigenvalue. If $N$
is large, we have $\lambda_i^N\ll \lambda_*^N$ for all $i\ge 3$, so
only the first two largest eigenvalues need to be taken into
account. This means that in order to evaluate $T^N$ we only need the
following eigenvectors: ${\bf v}_1=(1/n,1/n,\ldots, 1/n)^T$, ${\bf
v}_2$, $ {\bf w}_1=(n,0,0,\ldots, 0)$, and $ {\bf w}_2$.  The result
is:
\begin{equation}
\label{conv} 
Q_{11}=1-C\lambda_*^N,
\end{equation}
where $C=-\sum_{j=1}^n[{\bf v_2}]_j[{\bf w}_2]_1/n$. It follows,
therefore, that the convergence rate of the memoryless learner
algorithm can be estimated if we estimate $\lambda_*$ and $C$. It
turns out that once we understand $\lambda_*$ , we can also estimate
$C$.

Our results can be summarized as follows. For large $n$, the quantity
$C$ is bounded from above and below by some constants. From formulas
(\ref{muandlamb}) and (\ref{conv}) we can see that in order for the
learner to pick up the correct set with probability $1-\Delta$, we
need to have at least
\begin{equation}
\label{Nandmu} N_\Delta\sim|\log \Delta|/\mu_*
\end{equation}
sampling events (Theorem \ref{iweak} tells us that $\mu_* = o(1/n),$
and so we have the right to replace $\log(1-\mu_*)$ by $-\mu_*$). Using
the relationship between $\mu_*$ and the 
harmonic mean (\ref{Htomu}), and our results for $H_n$ from Theorem
\ref{iweak}, we obtain the following estimate:
\begin{equation}
\label{last} N_\Delta\sim |\log \Delta|h(n), 
\end{equation}
where $h(n)$ is $n \log n$ if the overlaps are uniformly distributed
(in other words, the entries $a_{1j}$ of the similarity matrix, as
random variables, are uniformly distributed in $[0, 1]$), and $h(n)$
is $n$ if the density of overlaps at $1$ goes to $0$. Estimate
(\ref{last}) should be understood in the sense that the right hand
side of (\ref{Nandmu}) converges in probability to the right hand side
of (\ref{last}). If the density grows at $1$ as
$(1-x)^\beta, -1 <\beta < 0$, then 
$$\lim_{x \rightarrow \infty}  \mathbf{P}\left(\frac{1}{x} < \frac{N_\Delta}{|\log
\Delta| n^{1/(1+\beta)}} < x\right) = 1.$$

\subsection{A better algorithm}
\label{fullmem}
Consider the following improvement on the previous learning algorithm:
the student keeps a list of the sets he has not rejected, and when the
time comes to switch, he picks uniformly among those sets
\emph{only}. It is clear that this algorithm (''learning with full
memory'') should perform better than the memoryless learner algorithm
described in the last section, but how much better?

Since the analysis is quite simple, we present it here. There are 
two questions which need to be answered (we always assume that 
the correct answer is the first set, $G_1$):

\begin{quote}
\textbf{Question 1.} Suppose the student has picked the set 
$G_i$, $i\neq 1$. What is the expected number of turns before 
he is forced to reject $G_i$ and jump to a different set? 
\end{quote}

\begin{quote}
\textbf{Question 2.} What is the probability that the student 
will change his mind exactly $k$ times before guessing the right 
answer? 
\end{quote}

We answer the second question first, by
\thm{lemma}{
The probability that the set $G_1$ is encountered on the $k$-th 
turn is independent of $k$ (and so equals $1/n$.)}
\begin{proof}
Suppose the student starts by picking a set $G_{i_1}$ at random, and
then keeps picking sets $G_{i_2}, G_{i_3}, \dots, G_{i_n}$, until
there are none left, and making sure never to repeat a set.  The
sequence $i_1, \dots, i_n$ is a permutation of the sequence $1, \dots,
n$, and it is clear (for reasons of symmetry) that every permutation
is equally likely. Since for any $k$, precisely $(n-1)!$ permutations
have $1$ in the $k$-th position, the lemma is proved.
\end{proof}
Question $1$ is also easily answered, by
\begin{lemma}
If $\nu_1(G_i) = a_i$, then the expected number of turns 
before switching is $1/(1-a_i)$.
\end{lemma}
\begin{proof}
Let $\mathcal{P}_k$ be the probability of switching on the $k$-th 
step or earlier. Then we have the equation: 
\begin{equation}
\mathcal{P}_{k+1} = \mathcal{P}_k + (1-\mathcal{P}_k) (1-a_i) = 
a_i \mathcal{P}_k  + (1-a_i). 
\end{equation}
Since $\mathcal{P}_0 = 0$, it is easy to check that $\mathcal{P}_j = 1
- a_i^j$. If $p_k$ is the probability of switching on the $k$-th
turn, then $p_k = a_i^{k-1} - a_i^k$, and the expected time of
switching is
\begin{equation}
\sum_{j=1}^\infty j (a_i^{j-1} - a_i^j) = \sum_{j=0}^\infty 
a_i^j = \frac{1}{1-a_i}, 
\end{equation}
the first equality being obtained by telescoping the sum.
\end{proof}
From the two lemmas, it follows that given the probabilities $a_2,
\dots, a_n$, the expected time taken by the improved learner is
\begin{equation*}
T= \frac{1}{n}\sum_{k=1}^{n-1}\frac{1}{\binom{n-1}{k}} 
\sum_{S_k}
\sum_{i\in S_k} \frac{1}{1-a_{i}}, 
\end{equation*}
where the middle summation is over all subsets $S_k$ of $2,\ldots, n$
which have size $k$. Since for any $i$, the number of subsets of $2, \dots, n$
of size $k$ containing $i$ equals $\binom{n-2}{k-1}$, the above
expression can be rewritten as
\begin{equation}
\begin{split}
T &= \frac{1}{n} \sum_{k=1}^{n-1} 
\frac{\binom{n-2}{k-1}}{\binom{n-1}{k}} 
\sum_{i=2}^{n}\frac{1}{1-a_i}\\ &= \sum_{i=2}^{n} \frac{1}{1-a_i} 
\sum_{k=1}^{n-1} \frac{k}{n(n-1)} =\frac{1}{2}\sum_{i=2}^n 
\frac{1}{1-a_i}=\frac{n-1}{2H_{n-1}},
\end{split}
\end{equation}
where $H_{n-1}$ is defined in (\ref{defharmean}) with
$x_i=1-a_i$. These computations can be easily adapted to solve the
following problem: suppose that we want to be $1-\Delta$ sure of
getting to the right answer. How many steps do we need? Notice that we
will need to take $(1-\Delta)n$ jumps, so the computation as above
gives us:
\begin{equation}
\begin{split}
N_\Delta &= \frac{1}{n} \sum_{k=1}^{(1-\Delta)n} 
\frac{\binom{n-2}{k-1}}{\binom{n-1}{k}} 
\sum_{i=2}^n\frac{1}{1-a_i}\\ &= \sum_{i=2}^n \frac{1}{1-a_i} 
\sum_{k=1}^{(1-\Delta)n} \frac{k}{n(n-1)}
\to \frac{n(1-\Delta)^2}{2H_{n}}. 
\end{split}
\end{equation}
Comparing this with equation (\ref{Nandmu}) and using estimate
(\ref{Htomu}), we notice that for every fixed $\Delta < 1$, this is
only a constant factor better than a memoryless learner. The constant
is a function of $\Delta$, and behaves as $\frac {\vert\log
\Delta\vert}{(1-\Delta)^2} \asymp \vert \log \Delta \vert,$ so goes to
infinity (albeit slowly) as $\Delta$ approaches $0$.

\end{document}